\documentclass[11pt]{amsart}

\usepackage{enumerate}

\numberwithin{equation}{section}

\newcommand{\goto}{\rightarrow}

\newcommand{\bg}{\bar {g}}

\newcommand{\cB}{\mathcal{B}}
\newcommand{\cF}{\mathcal{F}}
\newcommand{\cP}{\mathcal{P}}
\newcommand{\bN}{\mathbb{N}}
\newcommand{\bP}{\mathbb{P}}
\newcommand{\bR}{\mathbb{R}}
\newcommand{\bZ}{\mathbb{Z}}

\newcommand{\CZZ}{\hat {\mathbb{Z}}}

\newcommand{\XX}{\mathbf{X}}
\newcommand{\YY}{\mathbf{Y}}
\newcommand{\xx}{\mathbf{x}}
\newcommand{\yy}{\mathbf{y}}
\newcommand{\tx}{{\tilde{x}}}
\newcommand{\ty}{{\tilde{y}}}
\newcommand{\ta}{{\tilde{a}}}
\newcommand{\tb}{{\tilde{b}}}

\newcommand{\ee}{\mathbf{e}}

\newcommand{\al}{\alpha}
\newcommand{\be}{\beta}
\newcommand{\de}{\delta}
\newcommand{\De}{\Delta}

\theoremstyle{plain} \newtheorem{Theo}{Theorem}[section]
\theoremstyle{plain} \newtheorem{Lemma}[Theo]{Lemma}
\theoremstyle{plain} 
\theoremstyle{plain} \newtheorem{Prop}[Theo]{Proposition}
\theoremstyle{remark} \newtheorem{Rem}[Theo]{Remark}
\theoremstyle{definition} 
\theoremstyle{definition} 
\theoremstyle{plain}

\begin{document}

\centerline{\Large\bf A superprocess involving both branching and
coalescing}

\bigskip

\centerline{Xiaowen Zhou\footnote{\em E-mail address:
zhou@alcor.concordia.ca}}

\medskip

\centerline{Department of Mathematics and Statistics, Concordia
University,} \centerline{7141 Sherbrooke St. W.,  Montreal,  H4B
1R6, Canada}

\medskip

\centerline {\bf Abstract}

We consider a superprocess with  coalescing Brownian spatial motion.
We first prove a dual relationship between two systems of coalescing
Brownian motions. In consequence we can express the Laplace
functionals for the superprocess  in terms of coalescing Brownian
motions, which allows us to obtain some explicit results. We also
point out several connections between such a superprocess and the
Arratia flow. A more general model is discussed at the end of this
paper.

\medskip

\medskip



\section{Introduction}

In this paper we mainly consider the following branching-coalescing
particle system which can be described intuitively as follows. A
collection of particles with masses execute  coalescing Brownian
motions. In the mean while the masses for these particles evolve
according to independent  Feller's branching processes. Upon
coalescing those particles involved merge into one particle with
their respective masses added up.

The above-mentioned particle system can be described using a
measure-valued process $Z$. More precisely, the support of $Z_t$
represents the locations of those particles at time $t$, and the
measure $Z_t$ assigns to each supporting point stands for the mass
for the corresponding particle. This processes $Z$, referred as the
{\it superprocess with coalescing Brownian spatial motion} (SCSM),
was first introduced in \cite{DLZ04}. It arises as a scaling limit
of another measure-valued process, which was referred in
\cite{DLW01} as the {\it superprocess with dependent spatial motion}
(SDSM). As to SDSM, it arises as a high density limit of a critical
branching particle system in which the motion of each particle is
subjected to both an independent Brownian motion and a common white
noise applied to all the particles.  More precisely, the movement of
the $i^{\text{th}}$ particle is governed by equation
\[dx_i(t)=\sigma(x_i(t)) dB_i(t)+\int_\bR h(y-x_i(t))W(dy,dt),\]
where $(B_i)$ is a collection of independent Brownian motions which
is independent of the white noise $W$; see \cite{DLW01}. A similar
model was also studied in \cite{SkAd02}.

It was shown in Theorem 4.2 of \cite{DLZ04} that,
 after appropriate time-space scaling,  SDSM converges weakly
to SCSM. A functional dual for SCSM was  given in Theorem 3.4 of
\cite{DLZ04}. In addition,  using coalescing Brownian motions and
excursions for Feller's branching process, a construction of SDSM
was found in \cite{DLZ04},  an idea that initially came from
\cite{DL}. In this paper we always denote such a SCSM as $Z$.

One of the most interesting problems in the study of a
measure-valued process is to recover a certain dual relationship
concerning the measure-valued process. Such a dual relationship
often leads to the uniqueness of the measure-valued process; see
\cite{Per00} for some classical examples on {\it super Brownian
motion} and related processes. It is not hard to show the existence
of $Z$ as a high density limit of the branching-coalescing particle
system. The main goal of this paper is to propose a new way of
characterizing the measure-valued process $Z$ via duality, in which
the self duality for coalescing Brownian motions plays a key role.
To this end, we first prove a rather general duality on two
coalescing Brownian motions running in the opposite directions. We
derive this duality from an analogous, essentially combinatorial,
fact about coalescing simple random walk. With this duality we can
express certain Laplace functionals for $Z$ in terms of systems of
coalescing Brownian motions.

We could carry out some explicit computation thanks to the
above-mentioned  duality. In particular, we first show that,
starting with a possibly diffuse initial finite measure $Z_0$, $Z_t$
 collapses into a discrete measure with a finite support as soon as
$t>0$. Then we can identify $Z_t$ interchangeably with a finite
collection of spatially distributed particles with masses. When
there is such a particle at a fixed location, we  obtain the Laplace
transform of its mass. The total number of particles in $Z_t$
decreases in $t$ due to both branching and coalescing. When there is
only one particle left at time $t$, we  also recover the joint
distribution of its location and its mass. Eventually, all the
particles will die out. We  further find the distribution of the
location where the last particle disappears. Coincidentally, super
Brownian motion shares the same near extinction behavior.

Connections between superprocesses and stochastic flows have been
noticed before. In \cite{MaXi02} a superprocess was obtained from
the empirical measure of a coalescing flow. Arratia flow serves as a
fundamental example of coalescing flow. In this paper we point out
several connections between $Z$ and the Arratia flow. More
precisely, the support of $Z_t$ at a fixed time $t>0$ can be
identified with a Cox process whose intensity measure is determined
by the Arratia flow. A version of $Z_t$ can be constructed using the
Arratia flow. The general Laplace functional for $Z$ can also be
expressed in terms  of the Arratia flow.

Replacing the Feller's branching process by the square of Bessel
process, we discuss a more general model at the end of this paper.
The mass-dimension evolution of such a model can also be
characterized by coalescing Brownian motions.

The rest of this paper is arranged as follows.  As a preliminary, we
first state and prove a dual relationship on coalescing Brownian
motions in Section \ref{coal-dual}. In Section \ref{ex-uniq}, we
define the process $Z$  as a weak limit of the empirical measure for
the branching-coalescing particle system. Then we proceed to prove
the duality between $Z$ and coalescing Brownian motions. The
uniqueness of $Z$ follows from such a duality immediately. We
continue to study several properties of this process  in Section
\ref{prop}. We further discuss the connections between the Arratia
flow and $Z$ in Section \ref{flow}. At the end, we propose
 a more general model and establish its duality in Section 6.

\section{Coalescing Brownian motions and their duality}
\label{coal-dual}

An $m$-dimensional {\it coalescing Brownian motion} can be described
as follows. Consider a system of $m$ indexed particles with
locations in $\bR$ that evolve as follows. Each particle moves
according to an independent standard Brownian motion on $\bR$ until
two particles are at the same location. At this moment a {\em
coalescence} occurs and the particle of higher index starts to move
together with the particle of lower index. We say the particle with
higher index is {\em attached} to the particle with lower index,
which is still {\em free}. The particle system then continues its
evolution in the same fashion.  Note that indices are not essential
here, the collection of locations of the particles is Markovian in
its own right, but it will be convenient to think of the process as
taking values in $\bR^m$ rather than subsets of $\bR$ with at most
$m$ elements.  For definiteness, throughout this section we will
further assume that the particles are indexed in increasing order of
their initial positions: it it clear that the dynamics preserve this
ordering. Call the resulting Markov process $\XX = (X_1, \ldots,
X_m)$.

Write $1\{B\}(.)$ for the indicator function of a set $B$. The
distribution of $\XX(t)$ is uniquely specified by knowing for each
choice of $y_1 < y_2 < \ldots < y_n$ the joint probabilities of
which ``balls'' $X_1(t), X_2(t), \ldots, X_m(t)$ lie in which of the
``boxes'' $]y_1, y_2], ]y_2, y_3], \ldots, ]y_{n-1}, y_n]$. That is,
the distribution of $\XX(t)$ is determined by the joint distribution
of the indicators
\[
I_{ij}^\rightarrow(t,\yy) := 1\{X_i(t) \in ]y_{j}, y_{j+1}]\}
\]
for $1 \le i \le m$, $1 \le j \le n-1$ and $\yy=(y_1,\ldots,y_n)$.

Suppose now that $\YY := (Y_1, \ldots, Y_n)$ is another coalescing
Brownian motion.  The distribution of $\YY(t)$ is uniquely specified
by knowing for each choice of $x_1 < x_2 < \ldots < x_n$ the
distribution of the indicators
\[
I_{ij}^\leftarrow(t,\xx) := 1\{x_i \in ]Y_{j}(t), Y_{j+1}(t)]\}
\]
for $1 \le i \le m$, $1 \le j \le n-1$ and $\xx=(x_1,\ldots,x_m)$.

The next ``balls-in-boxes'' duality is crucial in characterizing the
distribution of the measure-valued process concerned in this paper.

\begin{Theo}\label{dual}
Suppose in the notation above that $\XX=(X_1,\ldots,X_m)$ is an
$m$-dimensional coalescing Brownian motion and
$\YY=(Y_1,\ldots,Y_n)$ is an $n$-dimensional coalescing Brownian
motion. Then for each $t \ge 0$ the joint distribution of the $m
\times (n-1)$-dimensional random array
$(I_{ij}^\rightarrow(t,\YY(0)))$ coincides with that of the $m
\times (n-1)$-dimensional random array
$(I_{ij}^\leftarrow(t,\XX(0)))$.
\end{Theo}

Theorem \ref{dual} generalizes Theorem 1.1 in \cite{XiZh05}. Some
other more elaborate dualities on coalescing-reflecting Brownian
systems can be found in \cite{ToWe97} and \cite{STW00}.

We first prove the counterpart of Theorem \ref{dual} for continuous
time {\it simple coalescing random walks}, which is  interesting in
its own right. Notice that $\XX$ is a coalescing Brownian motion if
and only if $X_i$ is a $(\cF^\XX_t)$-Brownian motion for each $1\leq
i\leq m$, and $(X_j-X_i)/\sqrt{2}$ is a $(\cF^\XX_t)$-Brownian
motion stopped at $0$, where $(\cF^\XX_t)$ denotes the filtration
generated by $\XX$. Then Theorem \ref{dual} follows from a straight
forward martingale argument proof of the convergence of scaled
random walk to Brownian motion.

A  $p$-simple random walk on $\bZ$ is a continuous time simple
random walk that makes jumps at unit rate, and when it makes a jump
from some site it jumps to the right neighbor with probability $p$
and to the left neighbor with probability $1-p$. An $m$-dimensional
{\it $p$-simple coalescing random walk} is defined in the same way
as the coalescing Brownian motion at the beginning of this section.
When $p=\frac{1}{2}$ we just call this particle system a simple
coalescing random walk.

Some notation is useful to keep track of the interactions among the
particles in the coalescing system. Let $\cP_m$ denote the set of
{\em interval partitions} of the totality of indices $\bN_m := \{1,
\ldots, m\}$. That is, an element $\pi$ of $\cP_m$ is a collection
$\pi = \{A_1(\pi), \ldots, A_h(\pi)\}$ of disjoint subsets of
$\bN_m$ such that $\bigcup_i A_i(\pi) = \bN_m$ and $a<b$ for all
$a\in A_i$, $b\in A_j$, $i<j$. The sets $A_1(\pi), \ldots A_h(\pi)$
consisting of consecutive indices are the {\em intervals} of the
partition $\pi$. The integer $h$ is the {\em length} of $\pi$ and is
denoted by $l(\pi)$. Equivalently, we can think of $\cP_m$ as a set
of equivalence relations on $\bN_m$ and write $i \sim_\pi j$ if $i$
and $j$ belong to the same interval of $\pi \in \cP_m$. Of course,
if $i \sim_\pi j$, then $i \sim_\pi k \sim_\pi j$ for all $i \le k
\le j$.

Given $\pi\in \cP_m$, define
$$\alpha_i(\pi) := \min A_i(\pi)$$
to be the left-hand end-point of the $i^{\mathrm{th}}$ interval
$A_i(\pi)$. Put
\[
\bZ^m_\pi := \{(x_1, \ldots, x_m) \in \bZ^m : x_1 \le \ldots \le x_m
\text{ and } x_i = x_j \text{ if } i \sim_\pi j\}
\]
and
\[
\begin{split}
\CZZ_\pi^m:=
\{(x_1, \ldots, x_m) \in \bZ^m : x_1 \le \ldots \le x_m \text{ and } x_i = x_j
\text{ if and only if } i \sim_\pi j\}.\\
\end{split}
\]
Note that $\bZ^m$ is the disjoint union of the sets $\CZZ^m_\pi$,
$\pi\in\cP_m$.

Write $\XX=(X_1,\ldots,X_m)$ for the coalescing random walk. If
$\XX(t) \in \CZZ^m_\pi$, then the free particles at time $t$ have
indices $\alpha_1(\pi),\ldots,\alpha_{l(\pi)}(\pi)$ and the
$i^{\mathrm{th}}$ particle at time $t$ is attached to the free
particle with index
\[
\min\{j : 1\le j \le m, \, j \sim_\pi i\} = \max\{\alpha_k(\pi) :
\alpha_k(\pi) \le i\}.
\]

In order to write down the generator of $\XX$, we require a final
piece of notation.
 Let  $\{\ee_{i}^k: 1\leq i\leq k\}$ be the
set of coordinate vectors in $\bZ^k$; that is, $\ee_{i}^k$ is the
vector that has the $i^{\mathrm{th}}$ coordinate $1$ and all the
other coordinates  $0$. For $\pi \in \cP_m$, define a map
$K_\pi:\bZ^m_\pi \rightarrow \bZ^{l(\pi)}$ by
\[
K_\pi(\xx)=K_\pi(x_1,\ldots,x_m):=\left(x_{\alpha_1(\pi)},\ldots,x_{\alpha_{l(\pi)}(\pi)}\right)
\]
Notice that $K_\pi$ is a bijection between $\bZ_\pi^m$ and $\{x \in
\bZ^{ l(\pi)} : x_1 \le x_2 \le \ldots \le x_{l(\pi)}\}$, and we
write $K_\pi^{-1}$ for the inverse of $K_\pi$. For brevity, we will
sometimes write $\xx_\pi$ for $K_\pi(\xx)$.

Write $B(\bZ^m)$ for the collection of all bounded functions on
$\bZ^m$. The generator $G$ of $\XX$ is the operator $G:
B(\bZ^m)\rightarrow B(\bZ^m)$ given by
\begin{equation*}
\begin{split}
G f (\xx) &:=p \sum_{i=1}^{l(\pi)} f\circ
K_\pi^{-1}(\xx_\pi+\ee_{i}^{l(\pi)}) +(1-p)
\sum_{i=1}^{l(\pi)}f\circ
K_\pi^{-1}(\xx_\pi-\ee_{i}^{l(\pi)}) \\
& \quad - l(\pi)f\circ K_\pi^{-1}(\xx_\pi),
\quad f\in B(\bZ^m), \; \xx\in \CZZ_\pi^m, \; \pi\in\cP_m.\\
\end{split}
\end{equation*}
This expression is well-defined, because if $\xx \in \CZZ_\pi^m$,
then $\xx_\pi$, $\xx_\pi+\ee_{i}^{l(\pi)}$ and
$\xx_\pi-\ee_{i}^{l(\pi)}$ are all in $\{x \in \bZ^{ l(\pi)} : x_1
\le x_2 \le \ldots \le x_{l(\pi)}\}$.

\noindent {\bf Note:} From now on we will suppress the dependence on
dimension and write $\ee_{i}^{l(\pi)}$ as $\ee_{i}$.

Write $\bZ':=\bZ+\frac{1}{2}=\{i+\frac{1}{2}:i\in \bZ\}$. An
$n$-dimensional $q$-simple coalescing random walk on ${\bZ'}^n$ and
its generator $H$ can be defined in the obvious way. Such a process,
with $q=1-p$, will serve as the process dual to the $p$-simple
coalescing random walk on $\bZ^m$ in the following way.

Fix $\xx \in \bZ^m$ with $x_1 \le \ldots \le x_m$ and $\yy \in
{\bZ'}^n$ with $y_1 \le \ldots \le y_n$. Put
\[
I_{ij}^\rightarrow(t,\yy) := 1\{X_i(t) \in ]y_{j}, y_{j+1}]\}
\]
and
\[
I_{ij}^\leftarrow(t,\xx) := 1\{x_i \in ]Y_{j}(t), Y_{j+1}(t)]\}
\]
for $1 \le i \le m$ and $1 \le j \le n-1$.

\begin{Lemma}\label{dual-crw}
Suppose in the notation above that $\XX=(X_1,\ldots,X_m)$ is an
$m$-dimensional $\bZ^m$-valued $p$-simple coalescing random walk and
$\YY=(Y_1,\ldots,Y_n)$ is an $n$-dimensional $\bZ'^n$-valued
$(1-p)$-simple coalescing random walk. Then for each $t \ge 0$ the
joint distribution of the $m \times (n-1)$-dimensional random array
$(I_{ij}^\rightarrow(t,\YY(0)))$ coincides with that of the $m
\times (n-1)$-dimensional random array
$(I_{ij}^\leftarrow(t,\XX(0)))$.
\end{Lemma}

\begin{proof}
For a function $g: \{0,1\}^{m(n-1)} \rightarrow \bR$, a vector $\xx
\in \bZ^m$ with $ x_1 \le \ldots \le  x_m$, and a vector $ \yy \in
{\bZ'}^n$ with $ y_1 \le \ldots \le  y_n$, set
\begin{equation*}
\bg( \xx;  \yy) :=g\left(1\{] y_1, y_2]\}(x_1), \ldots, 1\{]
y_{n-1},y_n]\}(x_1), \ldots,1\{]y_1, y_2]\}(x_m),
\ldots,1\{]y_{n-1}, y_n]\}(x_m)\right).
\end{equation*}
We may assume that $\XX$ and $\YY$ are defined on the same
probability space $(\Omega, \cF, \bP)$. We need to show that
\begin{equation}\label{ta1}
\bP[\bg(\XX(t);\YY(0))]=\bP[\bg(\XX(0);\YY(t))].
\end{equation}

For $\xx\in \bZ^m$, put $\bg_{\xx}(\cdot):=\bg( \xx;\cdot)$, and for
$\yy\in \bZ'^n$, put $\bg_{\yy}(\cdot):=\bg(\cdot; \yy)$. In order
to establish (\ref{ta1}), it suffices by a standard argument (cf.
Section 4.4 in \cite{EtKu86}) to show that
\begin{equation}
\label{genduality} G(\bg_{\yy})(\xx) = H(\bg_{\xx})(\yy)
\end{equation}
(recall that $G$ and $H$ are the generators of $\XX$ and $\YY$,
respectively).

Fix $\xx\in\CZZ_{\pi}^m$ and $\yy\in\CZZ_{\varpi}'^n$ for some
$\pi\in\cP_m$ and $\varpi\in\cP_n$. Put
\[
I^+:=\{i: 1\leq i\leq l(\pi), \; x_{\alpha_i(\pi)}+\frac{1}{2}
=y_{\alpha_j(\varpi)} \text{\, for some}\, 1\leq j\leq l(\varpi)\}
\]
and
\[
I^-:=\{i: 1\leq i\leq l(\pi), \;
x_{\alpha_i(\pi)}-\frac{1}{2}=y_{\alpha_j(\varpi)} \text{\, for
some}\, 1\leq j\leq l(\varpi)\}.\] Similarly, put
\[
J^-:=\{j: 1\leq j\leq l(\varpi), \;
y_{\alpha_j(\varpi)}-\frac{1}{2}=x_{\alpha_i(\pi)} \text{\, for
some}\, 1\leq i\leq l(\pi)\}
\]
and
\[J^+:=\{j: 1\leq j\leq l(\varpi), \; y_{\alpha_j(\varpi)}+\frac{1}{2}
=x_{\alpha_i(\pi)} \text{\, for some}\, 1\leq i\leq l(\pi)\}.\]

Recall that $x_{\alpha_1(\pi)} < \ldots < x_{\alpha_{l(\pi)}(\pi)}$
and $y_{\alpha_1(\varpi)} < \ldots <
y_{\alpha_{l(\varpi)}(\varpi)}$. Therefore, for each $i \in I^+$
there is a unique $j \in J^-$ such that
$x_{\alpha_i(\pi)}+\frac{1}{2}=y_{\alpha_j(\varpi)}$ and {\em vice
versa}. Fix such a pair $(i,j)$. Observe that
\[\xx':=\xx+\sum_{k\in A_i(\pi)}\ee_{k}^m=K^{-1}_\pi(\xx_\pi+\ee_i)  \]
and
\[\yy':=\yy-\sum_{k\in A_j(\varpi)}\ee_{k}^n=K^{-1}_\pi(\yy_\varpi-\ee_j).  \]
 We are going to verify that
\begin{equation}\label{ta2}
\left(1\{]y_{j'}, y_{j'+1}]\}(x'_{i'})\right)=\left(1\{]y'_{j'},
y'_{j'+1}]\}(x_{i'})\right)
\end{equation}
by considering all the possible scenarios.

Given any $i'\in A_i(\pi)$ we have:
\begin{itemize}
\item
for $j'=\alpha_j(\varpi)-1$,
\begin{equation*}
\begin{split}
&1\{]y_{j'},y_{j'+1}]\}(x'_{i'})=1\{]y_{j'},y_{j'+1}]\}(x_{i'}+1) \\
&\quad=0\\
&\quad=1\{]y_{j'},y_{j'+1}-1]\}(x_{i'})=1\{]y'_{j'},y'_{j'+1}]\}(x_{i'}),
\end{split}
\end{equation*}
\item
for $\alpha_j(\varpi)\leq j'<\max A_j(\varpi)$,
\begin{equation*}
\begin{split}
&1\{]y_{j'},y_{j'+1}]\}(x'_{i'})=1\{]y_{j'},y_{j'+1}]\}(x_{i'}+1)\\
&\quad=0\\
&\quad=1\{]y_{j'}-1,y_{j'+1}-1]\}(x_{i'})
=1\{]y'_{j'},y'_{j'+1}]\}(x_{i'}),
\end{split}
\end{equation*}
\item for $j'=\max A_j(\varpi)$,
\begin{equation*}
\begin{split}
&1\{]y_{j'},y_{j'+1}]\}(x'_{i'})=1\{]y_{j'},y_{j'+1}]\}(x_{i'}+1)\\
&\quad=1\\
&\quad=1\{]y_{j'}-1,y_{j'+1}]\}(x_{i'})=1\{]y'_{j'},y'_{j'+1}]\}(x_{i'}),
\end{split}
\end{equation*}
\item
and for $j'< \alpha_j(\varpi)-1$ or $j'> \max A_j(\varpi)$,
\begin{equation*}
\begin{split}
&1\{]y_{j'},y_{j'+1}]\}(x'_{i'})=1\{]y_{j'},y_{j'+1}]\}(x'_{i'})=1\{]y_{j'},y_{j'+1}]\}(x_{i'}+1)\\
&\quad=0\\
&\quad=1\{]y_{j'},y_{j'+1}]\}(x_{i'})=1\{]y'_{j'},y'_{j'+1}]\}(x_{i'}).
\end{split}
\end{equation*}
\end{itemize}
Moreover, given any $i'\not\in A_i(\pi)$, we have $x_{i'}\neq
x_{\alpha_i(\pi)}$. Hence
\begin{itemize}
\item for $j'=\alpha_j(\varpi)-1$,
\[1\{]y_{j'},y_{j'+1}]\}(x'_{i'})=1\{]y_{j'},y_{j'+1}]\}(x_{i'})=1\{]y_{j'},y_{j'+1}-1]\}(x_{i'})
=1\{]y'_{j'},y'_{j'+1}]\}(x_{i'}),\]
\item for $j'=\max A_j(\varpi)$,
\[1\{]y_{j'},y_{j'+1}]\}(x'_{i'})=1\{]y_{j'},y_{j'+1}]\}(x_{i'})=1\{]y_{j'}-1,y_{j'+1}]\}(x_{i'})
=1\{]y'_{j'},y'_{j'+1}]\}(x_{i'}),\]
\item for $\alpha_j(\varpi)\leq
 j'< \max A_j(\varpi)$,
\[1\{]y_{j'},y_{j'+1}]\}(x'_{i'})=1\{]y_{j'},y_{j'+1}]\}(x_{i'})=1\{]y_{j'}-1,y_{j'+1}-1]\}(x_{i'})
=1\{]y'_{j'},y'_{j'+1}]\}(x_{i'}),\] \item and for
$j'<\alpha_j(\varpi)-1$ or $j'>\max A_j(\varpi)$,
\[1\{]y_{j'},y_{j'+1}]\}(x'_{i'})=1\{]y_{j'},y_{j'+1}]\}(x_{i'})
=1\{]y'_{j'},y'_{j'+1}]\}(x_{i'}).\]
\end{itemize}
Combining the above observations yields (\ref{ta2}).

Therefore,
\[\bg_\yy \circ
K_\pi^{-1}(\xx_{\pi}+\ee_i)=\bg_\xx \circ
K_\pi^{-1}(\yy_{\varpi}-\ee_j).\] Furthermore, it is easy to see for
$i'\not\in I^+$ that
\[\bg_\yy \circ K_\pi^{-1}(\xx_{\pi}+\ee_{i'})
= \bg_\yy \circ K_\pi^{-1}(\xx_{\pi})\] and for $j'\not\in J^-$ that
\[\bg_\xx \circ K_\varpi^{-1}(\yy_{\varpi}-\ee_{j'})
=\bg_\xx \circ K_\varpi^{-1}(\yy_{\varpi}).\]

Similarly, for any $i\in I^-$ there exists a unique $j\in J^+$ such
that $x_{\alpha_i(\pi)}-\frac{1}{2}=y_{\alpha_j(\varpi)}$ and {\em
vice versa}.  For such a pair $(i,j)$ we have
\[\bg_\yy \circ K_\pi^{-1}(\xx_{\pi}-\ee_i)=\bg_\xx \circ
K_\pi^{-1}(\yy_{\varpi}+\ee_j).\] Furthermore, we see for $i'\not\in
I^-$ that
\[\bg_\yy \circ K_\pi^{-1}(\xx_{\pi}-\ee_{i'})= \bg_\yy \circ K_\pi^{-1}(\xx_{\pi}) \]
and for $j'\not\in J^+$ that
\[\bg_\xx \circ K_\varpi^{-1}(\yy_{\varpi}+\ee_{j'})=\bg_\xx \circ K_\varpi^{-1}(\yy_{\varpi}).\]

Lastly, note that
\[
\bg_\yy \circ K_\pi^{-1}(\xx_{\pi}) = \bg(\xx;\yy)= \bg_\xx \circ
K_\varpi^{-1}(\yy_{\varpi})
\]
and so
\begin{equation*}
\begin{split}
&G(\bg_{\yy})(\xx)-H(\bg_{\xx})(\yy)\\
&\quad=p \sum_{i=1}^{l(\pi)} \left(\bg_{\yy}\circ
K_\pi^{-1}(\xx_\pi+\ee_i)
-\bg_{\yy}\circ K_\pi^{-1}(\xx_\pi)\right)\\
&\qquad + (1-p) \sum_{i=1}^{l(\pi)}\left(\bg_{\yy}\circ
K_\pi^{-1}(\xx_\pi-\ee_i)
-\bg_{\yy}\circ K_\pi^{-1}(\xx_\pi)\right)\\
&\qquad -p \sum_{j=1}^{l(\varpi)} \left(\bg_{\xx}\circ
K_\varpi^{-1}(\yy_{\varpi}-\ee_i)-\bg_{\xx}\circ K_\varpi^{-1}(\yy_\varpi)\right)\\
&\qquad -(1-p) \sum_{j=1}^{l(\varpi)} \left(\bg_{\xx}\circ
K_\varpi^{-1}(\yy_\varpi+\ee_i)-\bg_{\xx}\circ K_\varpi^{-1}(\yy_\varpi)\right)\\
&\quad=p\sum_{i\in I^+}\bg_{\yy}\circ K_\pi^{-1}(\xx_\pi+\ee_i)
-p\sum_{j\in J^-}\bg_{\xx}\circ K_\varpi^{-1}(\yy_\varpi-\ee_j)\\
&\qquad +(1-p)\sum_{i\in I^-}\bg_{\yy}\circ
K_\pi^{-1}(\xx_\pi-\ee_i)
-(1-p)\sum_{j\in J^+}\bg_{\xx}\circ K_\varpi^{-1}(\yy_\varpi+\ee_j)\\
&\quad=0,\\
\end{split}
\end{equation*}
as required.
\end{proof}

\begin{Rem}
For discrete time coalescing simple random walks such a duality is
evident from  Fig. 7 in \cite{STW00}. But the duality seems to be
less apparent for continuous time coalescing simple random walk.
\end{Rem}

\section{Existence and uniqueness}
\label{ex-uniq}

A construction of $Z$ was given in \cite{DLZ04} using Feller's
branching excursions. In this paper we adopt a weak convergence
approach, which is commonly used in the study of measure-valued
processes.

Recall that a nonnegative valued process $\xi$ is a {\it Feller's
branching process} with initial value $x \geq 0$ if it is the unique
strong solution to the following stochastic differential equation
\[\xi(t)=x+\int_0^t\sqrt{\gamma \xi(s)}dB(s), \]
where $\gamma$ is a positive constant and $B$ is a one-dimensional
Brownian motion. $\xi(t)$ is a martingale. It has a Laplace
transform
\begin{equation}\label{lap-tran}
\bP\left[\exp\{-\lambda\xi(t)\}\right] =\exp\left\{-\frac{2\lambda
x}{2+\lambda\gamma t}\right\};
\end{equation}
its extinction probability is given by
\[\bP\{\xi(t)=0\}=\exp\left\{-\frac{2x}{\gamma t}\right\};\]
see Section II.1 and II.5 in \cite{Per00}.

Observe that independent Feller's branching processes are {\it
additive}; i.e. if $\xi$ and $\eta$ are two independent Feller's
branching processes (with the same parameter $\gamma$), then
$\xi+\eta$ is also a Feller's branching process. This fact will be
used repeatedly in our discussions.

Write $M_F(\bR)$ for the space of finite measures on $\bR$ equipped
with the topology of weak convergence. Given any finite measure
$Z_0$ on $\bR$, put $\bar{z}:=Z_0(\bR)$.  For any positive integer
$m$, let $(\xi_1^{(m)},\ldots,\xi_m^{(m)})$ be a collection of $m$
independent Feller's branching processes each with initial value
$\bar{z}/m$. Choose $(x_1,\ldots,x_n)$  to be i.i.d. samples from
distribution $Z_0/\bar{z}$. Let $(X_1^{(m)},\ldots,X_{m}^{(m)})$ be
an $m$-dimensional coalescing Brownian motion starting at
$(x_1,\ldots,x_m)$. Moreover, we always assume that
$(\xi_1^{(m)},\ldots,\xi_m^{(m)})$ and
$(X_1^{(m)},\ldots,X_m^{(m)})$ are independent. Let $\de_x$ denote
the point mass at $x\in\bR$.
 Then
\[Z_t^{(m)}:=\sum_{i=1}^m \xi_i^{(m)}(t)\delta_{X_i^{(m)}(t)}\]
defines a $M_F(\bR)$-valued process. From now on we will suppress
the dependence of $m$ in $\xi_i^{(m)}$ and $X_i^{(m)}$.

Recall that a collection of processes $\{Z^{\alpha}, \alpha\in I\}$
with sample paths in $D(M_F(\bR))$ is {\it C-relatively compact} if
it is relatively compact and all its weak limits are a.s.
continuous.  The proof of the next lemma is standard; see, e.g. the
proofs for Lemma 3.2 in \cite{XiZh05} and Proposition II.4.2 in
\cite{Per00}.

\begin{Lemma}\label{tight}
 $\{Z^{(m)}\}$ is C-relatively compact.
\end{Lemma}

\begin{proof}

We first check the {\it compact containment condition}. For any
$\epsilon
>0$ and $T>0$, choose a compact set $K_0\subset D(\bR)$ such that
$\bP\{X_1\in K_0^c\}< \epsilon^2 $. Let $K:=\{x_t: x\in K_0, t\leq
T\} $. Then $K$ is compact in $\bR$, and
\[\bP\{X_1(t)\in K^c, \exists t\leq T\}\leq \bP\{X_1\in
K_0^c\}<\epsilon^2 .\] Write
\[N_m:=\# I_K:=\#\{1\leq i\leq m: X_i(t)\in K^c, \exists t\leq T\},\]
where $\# I_K$ denotes the cardinality of the indices set $I_K$.
Conditioning on $N_m$, by the additivity for Feller's branching
processes, we see that $\sum_{i\in I_K}\xi_i $ is a Feller's
branching Process with initial value $N_m\bar{z}/m$. Then by Doob's
maximal inequality,
\[\bP\left\{\left.\sup_{0\leq t\leq T} \sum_{i\in I_K}\xi_i(t)> \epsilon\right|N_m \right\}
\leq \frac{N_m \bar{z}}{m\epsilon}.\]
 Therefore,
\begin{equation*}
\begin{split}
\bP\left\{\sup_{0\leq t\leq T} Z^{(m)}_t(K^c)>\epsilon\right\}
\leq\bP\left\{\sup_{0\leq t\leq T} \sum_{i\in
I_K}\xi_i(t)>\epsilon\right\} &\leq
\frac{\bP[N_m]\bar{z}}{m\epsilon}\leq \bar{z}\epsilon.
\end{split}
\end{equation*}

For any $ f\in C_b^2(\bR)$ put $ Z^{(m)}_t(f):=\int_{-\infty}^\infty
f(x) Z^{(m)}_t(dx)$. Now we are going to show that
$\{Z^{(m)}_.(f)\}$ is C-relatively compact in $D( \bR)$. By
It\^{o}'s formula, we have
\begin{equation*}
\begin{split}
Z_t^{(m)}(f)&=\sum_{i=1}^m \left[\frac{\bar{z}}{m}f(x_i)+\int_0^t
f(X_i(s))d\xi_i(s)+\int_0^t\xi_i(s)f^\prime(X_i(s))dX_i(s)\right.\\
&\quad\quad\quad\left.+\frac{1}{2}\int_0^t \xi_i(s)
f^{\prime\prime}(X_i(s))ds \right].
\end{split}
\end{equation*}

The additivity for $(\xi_i)$ gives
\[\bP\left[\sup_{0\leq s\leq t}\sum_{i=1}^m \xi_i(s)\right]<\infty
\,\,\, \text{and} \,\,\, \bP\left[\sup_{0\leq s\leq t}\sum_{i,j=1}^m
\xi_i(s)\xi_j(s)\right]<\infty, \, t>0, \] then $\sum_{i=1}^m
\int_0^t \xi_i(s) f^{\prime\prime}(X_i(s))ds$ is C-relatively
compact following from the Arzela-Ascoli theorem and Proposition
VI.3.26 of \cite{JaSh87}.

Note that
\[\left\langle\sum_{i=1}^m \int_0^.
\xi_i(s)f^\prime(X_i(s))dX_i(s)\right\rangle_t
=\sum_{i,j=1}^m\int_0^t
\xi_i(s)\xi_j(s)f^\prime(X_i(s))f^\prime(X_j(s))d\langle X_i,
X_j\rangle_s,\] where $\langle X_i, X_j\rangle_s=s-T_{ij}\wedge s$
and $T_{ij}:=\inf\{s\geq 0: X_i(s)=X_j(s)\}$. By Arzela-Ascoli
theorem again, $\{\langle\sum_{i=1}^m \int_0^.
\xi_i(s)f^\prime(X_i(s))dX_i(s)\rangle_.\}$ is C-relatively compact.
Theorem VI.4.13  and Proposition VI.3.26 in \cite{JaSh87} then imply
that the collection of martingales $\left\{\sum_{i=1}^m \int_0^.
\xi_i(s)f^\prime(X_i(s))dX_i(s)\right\}$ is C-relatively compact.

Similarly, $\left\{\sum_{i=1}^m \int_0^. f(X_i(s))d\xi_i(s)\right\}$
is also C-relatively compact. Moreover,
\[\frac{1}{m}\sum_{i=1}^m f(x_i)\goto Z_0(f) \text{\,\,a.s.}.\]
$\{Z^{(m)}(f)\}$ is thus C-relatively compact.  Consequently, by
Theorem II.4.1 in \cite{Per00}  we can conclude that $\{Z^{(m)}\}$
is C-relatively compact.

\end{proof}

Write $Z$ for the weak limit of $\{Z^{(m)}\}$.  The Laplace
functional of $Z$ can be obtained  from the duality in Theorem
\ref{dual}. As a result, its uniqueness is settled.

In the sequel we always write $(Y_1,\ldots,Y_{2n})$ for a coalescing
Brownian motion starting at $(y_1,\ldots,y_{2n})$ with
$y_1\leq\ldots\leq y_{2n} $.

\begin{Theo}\label{lap-fun}
$\{Z^{(m)}\}$ has a unique weak limit $Z$ in $C(M_F(\bR))$. Given
$a_j>0, j=1,\ldots, n$, for any $y_1\leq y_2\leq \ldots\leq y_{2n}$
and any $t>0$, we have
\begin{equation}\label{lap3}
\begin{split}
&\bP\left[\exp\left\{-\sum_{j=1}^n
a_j Z_t(]y_{2j-1},y_{2j}])\right\}\right]\\
&\quad=\bP\left[\exp\left\{-\int_{-\infty}^\infty
Z_0(dx)\frac{2\sum_{j=1}^n
a_j1{\{]Y_{2j-1}(t),Y_{2j}(t)]\}}(x)}{2+\gamma t\sum_{j=1}^n
a_j1\{]Y_{2j-1}(t),Y_{2j}(t)]\}(x)}\right\}\right].
\end{split}
\end{equation}

\end{Theo}

\begin{proof}
We first condition on $(\xi_1(t),\ldots,\xi_m(t))$. By Theorem
\ref{dual} we can show that
\begin{equation}\label{lap5}
\begin{split}
&\bP\left[\left.\exp\left\{-\sum_{j=1}^n \sum_{i=1}^m
a_j\xi_i(t)1\{]y_{2j-1},y_{2j}]\}(X_i(t))\right\} \right|
(\xi_i(t))\right]\\
&\quad=\bP\left[\left.\exp\left\{-\sum_{j=1}^n \sum_{i=1}^m
a_j\xi_i(t)1\{]Y_{2j-1}(t),Y_{2j}(t)]\}(x_i)\right\} \right|
(\xi_i(t))\right],
\end{split}
\end{equation}
where $(Y_1,\ldots,Y_{2n})$ is independent of $(X_i)$ and $(\xi_i)$.


Now take expectations on both sides of (\ref{lap5}) and then
condition on $(x_i)$ and $(Y_i(t))$. Since
$\xi_1(t),\ldots,\xi_m(t)$ are independent of each other, and they
are independent of $(x_i)$ and $(Y_i(t))$, it follows from
(\ref{lap-tran}) that
\begin{equation}\label{lap6}
\begin{split}
&\bP\left[\exp\left\{-\sum_{j=1}^n a_j Z_t^{(m)}(]y_{2j-1},y_{2j}])\right\}\right]\\
&\quad=\bP\left[\bP\left[\left.\exp\left\{-\sum_{j=1}^n \sum_{i=1}^m
a_j\xi_i(t)1\{]Y_{2j-1}(t),Y_{2j}(t)]\}(x_i)\right\} \right| (x_i),
(Y_i(t))\right]\right]\\
&\quad=\bP\left[\prod_{i=1}^m\exp\left\{-\frac{2\bar{z}\sum_{j=1}^n
a_j 1\{]Y_{2j-1}(t),Y_{2j}(t)]\}(x_i)}{m\{ 2+\gamma
t\sum_{j=1}^n a_j1\{]Y_{2j-1}(t),Y_{2j}(t)]\}(x_i)\}} \right\} \right]\\
&\quad=\bP\left[\left(\frac{1}{\bar{z}} \int_{-\infty}^\infty
Z_0(dx)\exp\left\{-\frac{2\bar{z}\sum_{j=1}^n
a_j1\{]Y_{2j-1}(t),Y_{2j}(t)]\}(x)}{m(2+\gamma t\sum_{j=1}^n
a_j1\{]Y_{2j-1}(t),Y_{2j}(t)]\}(x))} \right\}\right)^m \right].
\end{split}
\end{equation}

Let $Z$ be any weak limit of $\{Z^{(m)}\}$.
Let $m\goto\infty$ in (\ref{lap6}). Then
\begin{equation}
\begin{split}
&\lim_{m\goto\infty}
\bP\left[\exp\left\{-\sum_{j=1}^n a_j Z_t^{(m)}(]y_{2j-1},y_{2j}])\right\}\right]\\
&\quad=\lim_{m\goto\infty}\bP\left[\left(1-\int_{-\infty}^\infty
Z_0(dx)\frac{2\sum_{j=1}^n
a_j1\{]Y_{2j-1}(t),Y_{2j}(t)]\}(x)}{m(2+\gamma
t\sum_{j=1}^n a_j1\{]Y_{2j-1}(t),Y_{2j}(t)]\}(x))}\right)^m\right]\\
&\quad=\bP\left[\exp\left\{-\int_{-\infty}^\infty
Z_0(dx)\frac{2\sum_{j=1}^n
a_j1\{]Y_{2j-1}(t),Y_{2j}(t)]\}(x)}{2+\gamma t\sum_{j=1}^n
a_j1\{]Y_{2j-1}(t),Y_{2j}(t)]\}(x)}\right\}\right].
\end{split}
\end{equation}

We still need to make sure that
\begin{equation}\label{weak-con}
\begin{split}
\bP\left[\exp\left\{-\sum_{j=1}^n a_j
Z_t(]y_{2j-1},y_{2j}])\right\}\right] =\lim_{m\goto\infty}
\bP\left[\exp\left\{-\sum_{j=1}^n a_j Z_t^{(m)}(]y_{2j-1},y_{2j}])\right\}\right].\\
\end{split}
\end{equation}
To this end we can suppose that $y_1<y_2<\ldots<y_{2n}$. Then for
small enough $\epsilon>0$, similar to (\ref{lap6}) we have
\begin{equation}\label{conver}
\begin{split}
&\bP\left[\exp\left\{-\sum_{j=1}^n a_j
Z_t^{(m)}(]y_{2j-1}+\epsilon,y_{2j}-\epsilon])\right\}\right]
-\bP\left[\exp\left\{-\sum_{j=1}^n a_j
Z_t^{(m)}(]y_{2j-1}-\epsilon,y_{2j}+\epsilon])\right\}\right]\\
&\leq 1-\bP\left[\exp\left\{-\sum_{j=1}^n a_j
Z_t^{(m)}\left(]y_{2j-1}-\epsilon,y_{2j-1}+\epsilon]\cup
]y_{2j}-\epsilon,y_{2j}+\epsilon]\right)\right\}\right]\\
&=1-\bP\left[\left(\frac{1}{\bar{z}}\int_{-\infty}^\infty
Z_0(dx)\exp\left\{\frac{-2\bar{z} \sum_{j=1}^n
a_j1\{]Y'_{2j-1}(t),Y''_{2j-1}(t)]\cup
]Y'_{2j}(t),Y''_{2j}(t)]\}(x)}{m(2+\gamma t\sum_{j=1}^n
a_j1\{]Y'_{2j-1}(t),Y''_{2j-1}(t)]\cup
]Y'_{2j}(t),Y''_{2j}(t)]\}(x))}\right\}\right)^m\right]\\
&\leq 1-\bP\left\{\cap_{j=1}^{2n} \{Y'_j(t)=Y''_j(t)\}\right\},
\end{split}
\end{equation}
where $(Y'_1,Y''_1,\ldots,Y'_{2n},Y''_{2n}) $ is a coalescing
Brownian motion starting at
$(y_1-\epsilon,y_1+\epsilon,\ldots,y_{2n}-\epsilon,y_{2n}+\epsilon)$.
Clearly, (\ref{conver}) converges (uniformly in $m$) to $0$ as
$\epsilon\goto 0+$. So, (\ref{weak-con}) holds.

It is clear that  the distribution of $Z$ is uniquely determined by
(\ref{lap3}). So, $Z$ is the unique weak limit of $\{Z^{(m)}\}$.

\end{proof}

The moments of $Z$ can be obtained immediately from (\ref{lap3}).

\begin{Prop}
For any $y_1\leq y_2\leq\ldots\leq y_{2n}$ and $t>0$, we have
\[\bP\left[Z_t\left(\sum_{j=1}^n a_j1\{]y_{2j-1},y_{2j}]\}\right)\right]
=\bP\left[Z_0\left(\sum_{j=1}^n
a_j1\{]Y_{2j-1}(t),Y_{2j}(t)]\}\right)\right]\] and
\begin{equation*}
\begin{split}
&\bP\left[Z_t^2\left(\sum_{j=1}^n a_j1\{]y_{2j-1},y_{2j}]\}\right)\right]\\
&\quad=\bP\left[Z_0^2\left(\sum_{j=1}^n
a_j1\{]Y_{2j-1}(t),Y_{2j}(t)]\}\right)\right] +\bP\left[\gamma
tZ_0\left(\sum_{j=1}^n
a_j1\{]Y_{2j-1}(t),Y_{2j}(t)]\}\right)\right].\\
\end{split}
\end{equation*}
\end{Prop}

\begin{Rem}
 Another consequence of duality (\ref{lap3}) is
that $Z$ is a Markov process; see Theorem 3.3 in \cite{XiZh05} for a
proof on a similar model.
\end{Rem}

{\it Martingale problem} is often used to characterize a
superprocess.
 $Z$ is the solution to the martingale problem (see
\cite{DLZ04}): for any $\phi\in C^2(\bR)$,
\[M_t(\phi) = Z_t(\phi) - Z_0(\phi) - \frac{1}{\,2\,}
\int_0^t Z_s(\phi^{\prime\prime}) ds, \quad t\geq 0,\]
 is a continuous martingale relative to $({\cF}_t)_{t\geq 0}$ with
quadratic variation process
\[ \langle M(\phi)\rangle_t =
\gamma\int_0^t Z_s(\phi^2) ds + \int_0^t
ds\int_{\Delta}\phi^\prime(x)\phi^\prime(y) Z_s(dx)Z_s(dy), \]
where $\Delta = \{(x,x): x\in \bR\}$.

But a remarkable feature of such a martingale problem is that its
solution is {\it not} unique. For example, let $\xi_1$ and $\xi_2$
be two independent  branching processes each with initial value $1$.
Let $B_1$ and $B_2$ be two independent Brownian motions. Assume that
$(\xi_1, \xi_2)$ and $(B_1, B_2)$ are independent. Then
$Z_t':=\xi_1(t)\delta_{B_1(t)}+\xi_2(t)\delta_{B_1(t)}$ is another
solution to this martingale problem; also see \cite{XiZh05} for a
similar counter example.

\section{Some properties}
\label{prop}

Our first result in this section is a straight forward consequence
of Theorem \ref{lap-fun}.

\begin{Prop}\label{laplace}
 For any
$y_1\leq y_2\leq\ldots\leq y_{2n}$ and $t>0$, we have
\begin{equation}\label{lap1}
\begin{split}
&\bP\left[\exp\left\{-\lambda\sum_{j=1}^n Z_t(]y_{2j-1},y_{2j}])\right\}\right]\\
&\quad=\bP\left[\exp\left\{-\frac{2\lambda}{2+\lambda\gamma
t}\sum_{j=1}^n Z_0(]Y_{2j-1}(t),Y_{2j}(t)])\right\}\right],
\lambda>0.
\end{split}
\end{equation}
\end{Prop}

\begin{proof}
Observe that $\sum_{j=1}^n 1\{]Y_{2j-1}(t),Y_{2j}(t)]\}(x)$ is
either $0$ or $1$, then (\ref{lap1}) follows readily from
(\ref{lap3}).

\end{proof}

Proposition \ref{laplace} allows us to carry out some explicit
computations on $Z$. First, by letting $\lambda\goto\infty$ we can
easily see that
\[\bP[Z_t((-\infty,\infty))=0]=\exp\left\{-\frac{2\bar{z}}{\gamma t}\right\}.\]
We are going to further study the probability that $Z_t$ does not
charge on an arbitrary finite interval. For any $x,y, a$ and $b$,
put
\[\tx:=\frac{x-y}{\sqrt{2}}, \ty:=\frac{x+y}{\sqrt{2}},
\ta:=\frac{a+b}{\sqrt{2}} \text{\,\,\, and\,\,\,}
\tb:=\frac{b-a}{\sqrt{2}}.\]

\begin{Prop}\label{extinct}
Given $a<b$ and $t>0$, we have
\begin{equation}\label{extinct1}
\begin{split}
\bP\{Z_t(]a,b])=0\} &=\frac{1}{2\pi t}\int_{-\infty}^\infty
dx\int_0^\infty dy \exp\left\{-\frac{2Z_0(]\tx,\ty])}{t\gamma}
-\frac{(x-\ta)^2}{2t}\right\}\\
&\qquad\qquad \left(\exp\left\{-\frac{(y-\tb)^2}{2t}\right\}
-\exp\left\{-\frac{(y+\tb)^2}{2t}\right\}\right)\\
&\quad +\frac{2}{\sqrt{2\pi t}}\int_{\tb}^\infty dx
\exp\left\{-\frac{x^2}{2t}\right\}.
\end{split}
\end{equation}
\end{Prop}

\begin{proof}
Let $\lambda\goto\infty$ in (\ref{lap1}) we have
\begin{equation*}
\begin{split}
\bP\left\{ Z_t(]a,b])=0\right\}
&=\bP\left[\exp\left\{-\frac{2}{\gamma
t}Z_0(]Y_1(t),Y_2(t)])\right\}\right]\\
&=\bP\left[\exp\left\{-\frac{2}{\gamma
t}Z_0(]Y_1(t),Y_2(t)])\right\}; Y_1(t)\neq Y_2(t)\right]
+\bP\{Y_1(t)=Y_2(t)\},
\end{split}
\end{equation*}
where $(Y_1,Y_2)$ is a coalescing Brownian motion starting from
$(a,b)$.

To find the distribution of $(Y_1,Y_2)$, one could rotate the
coordinate system anti-clockwise by $\pi/4$. Under the new
coordinate system $(Y_1,Y_2)$ becomes a process $(Y_1', Y_2')$ such
that $Y_1'$ is a Brownian motion starting at  $\ta$, $Y_2'$ is a
Brownian motion starting at  $\tb$ and stopped at $0$, and $Y_1'$
and $Y_2'$ are independent. So (\ref{extinct1}) just follows from
the reflection principle for Brownian motion.

\end{proof}

Write $S_t$ for the support of $Z_t$.  Intuitively, starting with
particles of a total initial  mass $Z_0(\bR)$, as soon as $t>0$ the
particles near $-\infty$ and $\infty$ will die out due to branching.
$Z_t$ is then expected to be supported by a finite set  because of
coalescence. The next two results concern the cardinality of $S_t$.

\begin{Prop}\label{support}
Given $a<b$ and $t>0$, we have
\begin{equation}\label{support1}
\begin{split}
&\bP[\#S_t\cap ]a,b]]\\
& =\frac{b-a}{\sqrt{\pi t}}-\frac{1}{\sqrt{2}\pi t^2}\int_a^b
dz\int_{-\infty}^\infty dx\int_0^\infty
dyy\exp\left\{-\frac{2Z_0(]\tx,\ty
])}{t\gamma}-\frac{(x-\sqrt{2}z)^2+y^2}{2t}\right\}.
\end{split}
\end{equation}

\end{Prop}

\begin{proof}
It is easy to see from (\ref{extinct1}) that for any $z\in\bR$,
\begin{equation}\label{support2}
\begin{split}
&\bP\{Z_t(dz)\neq 0\}\\
&\quad=\frac{dz}{\sqrt{\pi t}}-\frac{dz}{\sqrt{2}\pi
t^2}\int_{-\infty}^\infty dx\int_0^\infty dy
y\exp\left\{-\frac{2Z_0(]\tx,\ty
])}{t\gamma}-\frac{(x-\sqrt{2}z)^2+y^2}{2t}\right\}.
\end{split}
\end{equation}
Then (\ref{support1}) is obtained by taking integrations on both
sides of (\ref{support2}) from $a$ to $b$.

\end{proof}

\begin{Prop}\label{fino}
With probability $1$, $\#S_t<\infty, \, \forall\, t>0$.
\end{Prop}

\begin{proof}
Given $s>0$,  we first claim that $\bP[\#S_s]<\infty$ if $Z_0$ has a
bounded support. Suppose that $Z_0(]-b,b])=1$ for some $b>0$. Then
by (\ref{support2}),
\begin{equation*}
\begin{split}
\bP[\#S_s] &=\frac{1}{\sqrt{2}\pi s^2}\int_{-\infty}^\infty
dz\int_{-\infty}^\infty dx\int_0^\infty dy
y\left(1-\exp\left\{-\frac{2Z_0(]\tx,\ty
])}{s\gamma}\right\}\right)\exp\left\{-\frac{(x-\sqrt{2}z)^2+y^2}{2s}\right\}\\
&=\frac{1}{\sqrt{2\pi s}s}\int_{-\infty}^\infty dx\int_0^\infty dy
y\left(1-\exp\left\{-\frac{2Z_0(]\tx,\ty
])}{s\gamma}\right\}\right)\exp\left\{-\frac{y^2}{2s}\right\}\\
&\leq \frac{1}{\sqrt{2\pi s}s}\int_0^\infty dy
\int_{-y-\sqrt{2}b}^{y+\sqrt{2}b}dx
y\exp\left\{-\frac{y^2}{2s}\right\}\\
&< \infty.
\end{split}
\end{equation*}
Our claim is proved.

Now given any integer $j$, let $\eta_j(s)$ be the Feller's branching
process with initial value $\eta_j(0):=Z_0(]j,j+1])$. Since
\[\sum_{j=-\infty}^\infty \bP\{\eta_j(s)\neq 0 \}=\sum_{j=-\infty}^\infty
\left(1-\exp\left\{-\frac{2\eta_j(0)}{\gamma s}\right\}\right)\leq
\sum_{j=-\infty}^\infty \frac{2\eta_j(0)}{\gamma
s}=\frac{2\bar{z}}{\gamma s},\] by Borel-Cantelli lemma we have
that, with probability $1$, $\eta_j(s)\neq 0$ for only finitely many
values of $j$.

Therefore, for any $t>0$, with probability $1$, $Z_{t/2}$ must have
a bounded support. The Markov property for $Z$, together with the
claim from the first part of the proof,  implies that $\#S_t<\infty$
a.s..

Finally, by the Markov property for $Z$ we conclude that
$\bP\{\#S_t<\infty, \forall \, t>0\}=1 $.
\end{proof}

By Proposition \ref{fino}, as soon as $t>0$, $S_t$ becomes a finite
set.   For any $z\in S_t$, we associate it with  a particle located
at $z$ with mass $Z_t(\{z\})$. We can thus identify $Z_t$
interchangeably with a collection of spatially distributed particles
with masses. As time goes on, the total number of particles
decreases either because two ``alive'' particles  coalesce into one
particle, or because each particle  disappears due to its branching.

Since $\#S_t<\infty$, a small neighborhood of $z$ contains at most
one particle in $Z_t$. When there is such a particle, we want to
find the distribution of its mass. Formally, we are looking for
\[\bP\left[\exp\{-\lambda Z_t(\{z\}) \};Z_t(\{z\})>0\right].\]

\begin{Prop}
For any $z\in\bR$ and $t>0$, we have
\begin{equation}\label{joeq}
\begin{split}
&\bP\left[\exp\{-\lambda Z_t(dz) \};Z_t(dz)>0\right]\\
&=\frac{dz}{\sqrt{2}\pi t^2}\int_{-\infty}^\infty dx\int_0^\infty dy
y\left(\exp\left\{-\frac{2\lambda Z_0(]\tx,\ty])}{2+\lambda\gamma
t}\right\}-\exp\left\{-\frac{2Z_0(]\tx,\ty])}{t\gamma}\right\}\right)\\
&\qquad\quad \exp\left\{-\frac{(x-\sqrt{2}z)^2+y^2}{2t}\right\}.
\end{split}
\end{equation}
\end{Prop}

\begin{proof}

We fix $(\xi_i(t))$ first. Apply Theorem \ref{dual} to
\[\exp\left\{-\lambda\sum_{i=1}^m\xi_i(t)1\{]a,b]\}(X_i(t))\right\}
1\left\{ \sum_{i=1}^m\xi_i(t)1\{]a,b]\}(X_i(t))>0\right\}.\] Then
condition on $ (Y_1(t),Y_2(t))$ and take an expectation with respect
to $(\xi_i(t))$. Similar to the proof for Theorem \ref{lap-fun} we
have that
\begin{equation}
\begin{split}
&\bP\left[\exp\left\{-\lambda\sum_{i=1}^m\xi_i(t)1\{]a,b]\}(X_i(t))\right\};
\sum_{i=1}^m\xi_i(t)1\{]a,b]\}(X_i(t))>0\right]\\
&\quad=\bP\left[\exp\left\{-\lambda\sum_{i=1}^m\xi_i(t)1\{]Y_1(t),Y_2(t)]\}(x_i)\right\};
\sum_{i=1}^m\xi_i(t)1\{]Y_1(t),Y_2(t)]\}(x_i)>0\right]\\
&\quad=\bP\left[\exp\left\{-\lambda\sum_{i=1}^m\xi_i(t)1\{]Y_1(t),Y_2(t)]\}(x_i)\right\}\right]
-\bP\left[\sum_{i=1}^m\xi_i(t)1\{]Y_1(t),Y_2(t)]\}(x_i)=0\right]\\
&\quad=\bP\left[\exp\left\{-\frac{2\lambda
Z_0^{(m)}(]Y_1(t),Y_2(t)])}{2+\lambda\gamma t}\right\}\right]-
\bP\left[\exp\left\{-\frac{2Z_0^{(m)}(]Y_1(t),Y_2(t)])}{\gamma
t}\right\}\right].
\end{split}
\end{equation}
Therefore,
\begin{equation}
\begin{split}
&\bP\left[\exp\{-\lambda Z_t(]a,b]) \};Z_t(]a,b])>0\right]\\
&\quad=\bP\left[\exp\left\{-\frac{2\lambda
Z_0(]Y_1(t),Y_2(t)])}{2+\lambda\gamma t}\right\}\right]-
\bP\left[\exp\left\{-\frac{2Z_0(]Y_1(t),Y_2(t)])}{\gamma
t}\right\}\right]\\
&\quad=\frac{1}{2\pi t}\int_{-\infty}^\infty dx\int_0^\infty dy
\left(\exp\left\{-\frac{2\lambda Z_0(]\tx,\ty])}{2+\lambda\gamma
t}\right\}-\exp\left\{-\frac{2Z_0(]\tx,\ty])}{t\gamma}\right\}\right)\\
&\qquad\quad \exp\left\{-\frac{(x-\ta)^2}{2t}\right\}
\left(\exp\left\{-\frac{(y-\tb)^2}{2t}\right\}
-\exp\left\{-\frac{(y+\tb)^2}{2t}\right\}\right).
\end{split}
\end{equation}
So, (\ref{joeq}) is obtained by letting $b\goto a+$.

\end{proof}

At a fixed time $t>0$, with a positive probability there can be only
one particle (with a positive mass) left. When this happens, we are
interested in the joint distribution of the  mass and the location
of that particle. More precisely, we want to find
\[\bP\left[\exp\{-\lambda Z_t(\bR)\}; Z_t(\bR)\neq 0,
S_t\subset dz\right].\]

\begin{Prop}\label{joint}
For any $z\in\bR$ and $t>0$, we have
\begin{equation}\label{joint3}
\begin{split}
&\bP\left[\exp\{-\lambda Z_t(\bR)\}; Z_t(\bR)\neq 0,
S_t\subset dz\right]\\
&=\frac{dz}{\sqrt{2}\pi t^2}\int_{-\infty}^\infty dx\int_0^\infty dy
y\exp\left\{-\frac{2\lambda Z_0(]\tx,\ty])}{2+\lambda\gamma
t}-\frac{2Z_0(]\tx,\ty]^c)}{\gamma t}-\frac{(x-\sqrt{2}z)^2+y^2}{2t}\right\}\\
&\quad-\frac{dz}{\sqrt{\pi t}}\exp\left\{-\frac{2\bar{z}}{\gamma
t}\right\}.
\end{split}
\end{equation}
\end{Prop}

\begin{proof}

Put \[B:=\left\{\sum_{i=1}^m
\xi_i(t)1\{]Y_1(t),Y_2(t)]^c\}(x_i)=0\right\}\,\,\,\text{ for} \,\,
x_i:=X_i(0).\] It follows from Theorem \ref{dual} that
\begin{equation}\label{joint2}
\begin{split}
&\bP\left[\exp\left\{-\lambda\sum_{i=1}^m\xi_i(t)\right\};
\sum_{i=1}^m \xi_i(t)1\{]a,b]\}(X_i(t))>0, \sum_{i=1}^m \xi_i(t)1\{]a,b]^c\}(X_i(t))=0 \right]\\
&= \bP\left[\exp\left\{-\lambda\sum_{i=1}^m\xi_i(t) \right\};
\sum_{i=1}^m \xi_i(t)1\{]Y_1(t),Y_2(t)]\}(x_i)>0, B\right]\\
&=\bP\left[\exp\left\{-\lambda\sum_{i=1}^m\xi_i(t)1\{]Y_1(t),Y_2(t)]\}(x_i)
\right\}; \sum_{i=1}^m \xi_i(t)1\{]Y_1(t),Y_2(t)]\}(x_i)>0, B
\right]\\
&=\bP\left[\exp\left\{-\lambda\sum_{i=1}^m\xi_i(t)1\{]Y_1(t),Y_2(t)]\}(x_i)
\right\}; B\right]-\bP\left[\sum_{i=1}^m\xi_i(t)1\{]Y_1(t),Y_2(t)]\}(x_i)=0; B\right]\\
&=\bP\left[\exp\left\{-\frac{2\lambda \bar{z}\sum_{i=1}^m
1\{]Y_1(t),Y_2(t)]\}(x_i)}{m(2+\lambda\gamma t)}\right\}
\exp\left\{-\frac{2\bar{z}\sum_{i=1}^m
1\{]Y_1(t),Y_2(t)]^c\}(x_i)}{m\gamma
t}\right\}\right]\\
&\quad-\bP\left\{\sum_{i=1}^m\xi_i(t)=0\right\},
\end{split}
\end{equation}
where in obtaining the last equation we have used the fact that,
given $(Y_1(t),Y_2(t))$ and $(x_i)$,
\[\exp\left\{-\lambda\sum_{i=1}^m\xi_i(t)1\{]Y_1(t),Y_2(t)]\}(x_i)\right\}
\text{ \,\,\, and event\,\,\,\,}
\left\{\sum_{i=1}^m\xi_i(t)1\{]Y_1(t),Y_2(t)]^c\}(x_i)=0\right\}\]
are independent.

Now letting $m\goto\infty$ in (\ref{joint2}) we have
\begin{equation}\label{range}
\begin{split}
&\bP\left[\exp\{-\lambda Z_t(\bR)\}; Z_t(\bR)\neq 0,
S_t\subset (a,b)\right]\\
&=\bP\left[\exp\left\{-\frac{2\lambda
Z_0(]Y_1(t),Y_2(t)])}{2+\lambda\gamma t}
-\frac{2Z_0(]Y_1(t),Y_2(t)]^c)}{\gamma t}\right\}\right]-\exp\left\{-\frac{2\bar{z}}{\gamma t}\right\}\\
&=\frac{1}{2\pi t}\int_{-\infty}^\infty dx\int_0^\infty dy
\exp\left\{-\frac{2\lambda Z_0(]\tx,\ty])}{2+\lambda\gamma
t}-\frac{2Z_0(]\tx,\ty]^c)}{\gamma t}\right\}
\exp\left\{-\frac{(x-\ta)^2}{2t}\right\}\\
&\qquad\qquad \left(\exp\left\{-\frac{(y-\tb)^2}{2t}\right\}
-\exp\left\{-\frac{(y+\tb)^2}{2t}\right\}\right)\\
&\quad +\frac{2}{\sqrt{2\pi t}}\exp\left\{-\frac{2}{\gamma
t}\right\}\int_{\tb}^\infty dx
\exp\left\{-\frac{x^2}{2t}\right\}-\exp\left\{-\frac{2\bar{z}}{\gamma
t}\right\}.
\end{split}
\end{equation}
Finally, (\ref{joint3}) is obtained by letting $b\goto a+$ in
(\ref{range}).

\end{proof}

\begin{Rem}
Let $\lambda=0$ in (\ref{range}). We then obtain a result on the
range of $S_t$.
\end{Rem}

 The total number of particles in $Z$ will decrease one by one. Put
\[\tau:=\inf\{s\geq 0:\# S_s=1\}.\]
Then $\tau<\infty$ is the first time when there is exactly one
particle left. The distribution of $\tau$ is given in the following
Proposition.

\begin{Prop}
\begin{equation}\label{one-left}
\begin{split}
\bP\{\tau<t\} &=\int_{-\infty}^\infty dz\left(\frac{1}{\sqrt{2}\pi
t^2}\int_{-\infty}^\infty dx\int_0^\infty dy
y\exp\left\{-\frac{2Z_0(]\tx,\ty]^c)}
{\gamma t}-\frac{(x-\sqrt{2}z)^2+y^2}{2t}\right\}\right.\\
&\qquad\left.-\frac{1}{\sqrt{\pi
t}}\exp\left\{-\frac{2\bar{z}}{\gamma
t}\right\}\right)+1-\exp\left\{\frac{2\bar{z}}{\gamma t}\right\}.
\end{split}
\end{equation}
\end{Prop}

\begin{proof}
Observe that
\begin{equation*}
\begin{split}
\bP\{\tau<t\} &=\int_{-\infty}^\infty\bP\left\{Z_t(\bR)\neq 0,
S_t\subset dz\right\}
+\bP\{Z_t(\bR)= 0\},\\
\end{split}
\end{equation*}
then (\ref{one-left}) follows from Proposition \ref{joint}.

\end{proof}

Let
\[T:=\inf\{t\geq 0: Z_t(\bR)= 0\}.\]
$T$ is the time when all the particles disappear. Its distribution
can be found easily.
\[\bP\{T\leq t\}=\bP\{Z_t(\bR)=0\}=\exp\left\{-\frac{2\bar{z}}{\gamma t}\right\}. \]
Let $F$ denote the location of the last particle immediately before
extinction, i.e. $\{F\}=S_{T-} $. We could recover the explicit
distribution for $F$.

\begin{Prop}
$F$ has the same distribution as $X_{T}$, where $X$ is a Brownian
motion  with initial distribution $Z_0/\hat{z}$, and $X$ and $T$ are
independent.
\end{Prop}

\begin{proof}
First assume that \[Z_0=\sum_{i=1}^m a_i\delta_{x_i}\] with $a_i>0$
and $\sum_{i=1}^m a_i=a$. Then
\[Z_t:=\sum_{i=1}^m\xi_i(t)\delta_{X_i(t)},\] where
$\xi_i(0)=a_i>0$ and $(X_1,\ldots,X_m)$ is a coalescing Brownian
motion starting at $(x_1,\ldots,x_m)$.

Write $T_i:=\inf\{t\geq 0:\xi_i(t)=0\}, i=1,\ldots,m$. Then
$T:=\max_{1\leq i\leq m}T_i$. Therefore,
\[F=\sum_{i=1}^m X_i(T_i-)1\{T=T_i\}=\sum_{i=1}^m X_i(T)1\{T=T_i\}.\]

Our first observation is that
\[\bP\{T_i\leq t\}=\bP\{\xi_i(t)=0\}=\exp\left\{-\frac{2a_i}{\gamma t}\right\}.\]
Then $\bP\{T=T_i\}=a_i/a$, and $F=X_i(T)$ with probability $a_i/a$.
Our second observation is that conditional on $\{T=T_i\}$, the
distribution for $T$ is the same  as its unconditional distribution.
So, $F$ has the same distribution as the random variable obtained by
running a Brownian motion $X$ with initial distribution
$\bP\{X(0)=x_i\}=a_i/a, i=1,\ldots,m$, and stopping it independently
at time $T$. As a result, $F$ has the desired distribution.

By conditioning on $Z_\epsilon$ and letting $\epsilon\goto 0+$, the
conclusion in the proposition also follows for any general initial
measure $Z_0$.

\end{proof}


\begin{Rem}
This near extinction behavior is the same as  the super Brownian
motion (see Theorem 1 in \cite{Tri92}).
\end{Rem}

\section{Connections with the Arratia flow}
\label{flow}

Arratia flow is a stochastic flow which describes the evolution of a
continuous family of coalescing Brownian motions on $\bR$. We refer
to \cite{Ar79} for a detailed account and \cite{Da} for a survey on
stochastic flows. By definition, the Arratia flow  is a collection
$\{\phi(s,t,x):0\leq s\leq t, x\in \bR\}$ of random variables such
that

\begin{itemize}
\item the random map $(s,t,x) \mapsto \phi(s,t,x)$ is jointly
measurable, \item for each $s$ and $x$, the map $t \mapsto
\phi(s,t,x)$, $t \ge s$, is continuous,
\item for each $s$ and $t$
with $s \le t$, the map $x \mapsto \phi(s,t,x)$ is non-decreasing
and right-continuous, \item for $s \le t \le u$, $\phi(t,u,\cdot)
\circ \phi(s,t,\cdot) = \phi(s,u,\cdot)$,
\item for $u>0$,
$(s,t,x) \mapsto \phi(s+u, t+u, x)$ has the same distribution as
$\phi$, \item for $x_1< \ldots< x_m$ the process $(\phi(0,t,x_1),
\ldots, \phi(0,t,x_m))_{t \ge 0}$ has the same distribution as a
coalescing Brownian motion starting at $(x_1, \ldots, x_m)$.
\end{itemize}

Fix $t>0$, it is known that  $\{\phi(0,t,x): x\in\bR \}$, the image
of $\bR$ under map $\phi(0,t,.)$, is a discrete set (see
\cite{Ar79}). Let $\ldots<x^*_{-1}<x^*_0<x^*_1<\ldots$ be a sequence
of random variables such that
\begin{equation}\label{image}
\{\phi(0,t,x): x\in\bR \}=\{\ldots, x^*_{-1}, x^*_0, x^*_1,\ldots\}.
\end{equation}
Since Brownian motion has continuous sample paths, the Arratia flow
is order-preserving; i.e. $\phi(0,t,x_1)\leq\ldots\leq
\phi(0,t,x_m)$ whenever $x_1\leq,\ldots,\leq x_m$. Set
\[\Pi_i:=\sup\{x:\phi(0,t,x)=x^*_i
\}.\] Write $\phi^{-1}(0,t,x)$ for the pre-image of $x$ under
$\phi(0,t,.)$ Then $(\Pi_i)$ determines a partition on $\bR$ such
that $\phi^{-1}(0,t,x^*_i)=[\Pi_{i-1},\Pi_i[$.

Not surprisingly, the Arratia flow is closely connected to the
process $Z$ studied in the previous sections. We first consider its
support $S_t$. Since $S_t$ is a discrete set, we can identify it
with a simple point process by placing a unit mass on each point of
$S_t$. For any $y_1\leq y_2\leq\ldots\leq y_{2n}$, by Proposition
\ref{laplace},
\begin{equation}\label{occu}
\bP\left\{Z_t(\cup_{j=1}^n ]y_{2j-1},y_{2j}])=0\right\}
=\bP\left[\exp\left\{-\frac{2}{\gamma
t}Z_0(\cup_{j=1}^n]Y_{2j-1}(t),Y_{2j}(t)])\right\}\right].
\end{equation}
We thus get the following characterization of the {\it avoidance
function} for $S_t$.
\begin{equation}\label{supt}
\bP\left\{S_t\cap\cup_{j=1}^n ]y_{2j-1},y_{2j}]=\emptyset\right\}
=\bP\left[\exp\left\{-\frac{2}{\gamma
t}Z_0(\cup_{j=1}^n]Y_{2j-1}(t),Y_{2j}(t)])\right\}\right].
\end{equation}
Consequently the distribution of $S_t$ is uniquely determined by
(\ref{supt}); see Theorem 3.3 in \cite{Kal76}.

(\ref{supt})  suggests a connection between $S_t$ and the Arratia
flow. Let $M_t(dy)$ be a random measure on $\bR$ such that
\[M_t\left(\cup_{j=1}^n ]y_{2j-1},y_{2j}]\right)
=\frac{2}{\gamma t}\sum_{j=1}^n Z_0\left(]\phi(0,t,y_{2j-1}),
\phi(0,t,y_{2j})]\right), \,\,\, y_1\leq y_2\leq\ldots\leq y_{2n}.
\]
Then $S_t$ can be identified with a Cox process with a finite random
intensity measure $M_t$.

(\ref{occu}) also leads to a result on the {\it occupation time} for
$Z$. For any Borel set $B$ in $\bR$,
\[\int_0^tds\bP\{ Z_s(B)=0\}
=\int_0^tds\bP\left[\exp\left\{-M_t(B)\right\}\right].\]

A particle representation for $Z_t$ is available by using the image
of the Arratia flow as a skeleton. Given $(x^*_i)$ as in
(\ref{image}), let $(\ldots,\xi_{-1},\xi_0,\xi_1,\ldots)$ be
independent non-negative random variables such that
\[\bP\left[\left.\exp\{-\lambda\xi_i\}\right|(x^*_i)\right]
=\exp\left\{-\frac{2\lambda
Z_0(\phi^{-1}(0,t,x^*_i))}{2+\lambda\gamma t}\right\}.\] Then
\begin{equation}\label{particle-con}
 Z_t\overset{D}{=}\sum_{i=-\infty}^\infty \xi_i
\delta_{x^*_i}.
\end{equation}

To see this, define
\[Z^{(m)}_t:=\sum_{i=-m2^m}^{m2^m}\xi^{(m)}_i(t)\delta_{\phi\left(0,t,i/2^m\right)}, \]
where $\left(\xi^{(m)}_i\right)_{i=-m2^m}^{m2^m}$ is a sequence of
independent Feller's branching processes with initial values
$\left(Z_0([(i-1)/2^m,i/2^m[)\right)_{i=-m2^m}^{m2^m}$, and in
addition, $\left(\xi^{(m)}_i\right)$ is independent of
$\{\phi(s,t,x)\}$.

For any $a_j\geq 0, j=1,\ldots,n$, and $y_1\leq\ldots\leq y_{2n}$,
by the same argument as in the proof for Theorem \ref{lap-fun}, we
have
\begin{equation*}
\begin{split}
&\lim_{m\goto\infty}\bP\left[\exp\left\{-\sum_{j=1}^n a_j Z_t^{(m)}(]y_{2j-1},y_{2j}])\right\}\right]\\
&\quad=\lim_{m\goto\infty}\bP\left[\prod_{i=-m2^m}^{m2^m}
\exp\left\{-\frac{2Z_0([(i-1)/{2^m},i/{2^m}[) \sum_{j=1}^n
a_j 1\{]Y_{2j-1}(t),Y_{2j}(t)]\}(i/{2^m})}{2+\gamma
t\sum_{j=1}^n a_j 1\{]Y_{2j-1}(t),Y_{2j}(t)]\}(i/2^m)} \right\} \right]\\
&\quad=\lim_{m\goto\infty}\bP\left[\exp\left\{-\sum_{i=-m2^m}^{m2^m}\frac{2Z_0([(i-1)/2^m,i/2^m[)
\sum_{j=1}^n a_j 1\{]Y_{2j-1}(t),Y_{2j}(t)]\}(i/2^m)}{2+\gamma
t\sum_{j=1}^n a_j 1\{]Y_{2j-1}(t),Y_{2j}(t)]\}(i/2^m)}
\right\} \right]\\
&\quad=\bP\left[\exp\left\{-\int_{-\infty}^\infty
Z_0(dx)\frac{2\sum_{j=1}^n a_j
1\{]Y_{2j-1}(t),Y_{2j}(t)]\}(x)}{2+\gamma t\sum_{j=1}^n a_j
1\{]Y_{2j-1}(t),Y_{2j}(t)]\}(x)}\right\}\right]\\
&\quad=\bP\left[\exp\left\{-\sum_{j=1}^n a_j
Z_t(]y_{2j-1},y_{2j}])\right\}\right].
\end{split}
\end{equation*}
Therefore,
\[Z^{(m)}_t\overset{D}{\goto} Z_t.\]

Further, by the definition of $(x^*_i)$ and the additive property
for Feller's branching processes we obtain that
\[Z^{(m)}_t=\sum_{i=-\infty}^\infty\sum_{\Pi_{i-1}\leq {j}/{2^m}<\Pi_i}\xi^{(m)}_j(t)\delta_{x^*_i}
\overset{D}{\goto}\sum_{i=-\infty}^\infty \xi_i \delta_{x^*_i}.\]
Putting these together gives (\ref{particle-con}).

This interplay is remarkable. On one hand, $Z$ can be constructed
using the Arratia flow; on the other hand, $Z$ tells us how an
initial measure $Z_0$ is transported over time under both the
Arratia flow and the branching.

 The Laplace functional for
$Z_t$ can also be expressed in terms of $(x^*_i)$ and $(\Pi_i)$.
Given any nonnegative bounded continuous function $f$, for
$y_j=j/2^n$, Theorem \ref{lap-fun} yields
\begin{equation}\label{flow1}
\begin{split}
&\bP\left[\exp\left\{-\sum_{j=-n2^n}^{n2^n} f(y_j) Z_t(]y_{j-1},y_j])\right\}\right]\\
&\quad=\bP\left[\exp\left\{-\int_{-\infty}^\infty
Z_0(dx)\frac{2\sum_{j=-n2^n}^{n2^n}f(y_j) 1\{]\phi(0,t,y_{j-1}),
\phi (0,t,y_j)]\}(x)}{2+\gamma t\sum_{j=-n2^n}^{n2^n}
f(y_j)1\{]\phi(0,t,y_{j-1}), \phi (0,t,y_j)]\}(x)}\right\}\right].
\end{split}
\end{equation}
Again, let $m\goto\infty$. It follows that
\begin{equation*}
\begin{split}
&\bP\left[\exp\left\{-\lambda \int_{-\infty}^\infty
f(x)Z_t(dx)\right\}\right]\\
&\quad=\bP\left[\exp\left\{-\int_{-\infty}^\infty
Z_0(dx)\frac{2\lambda\sum_{i=-\infty}^\infty
f(\Pi_i)1\{]x^*_i,x^*_{i+1}]\}(x)}{2+\lambda\gamma t
\sum_{i=-\infty}^\infty
f(\Pi_i)1\{]x^*_i,x^*_{i+1}]\}(x)}\right\}\right].
\end{split}
\end{equation*}

\section{A more general model}
\label{general}

Evans observes that what is really at work in the proof for Theorem
\ref{lap-fun} is the additivity for the Feller's branching
processes. He then suggested that we could use the square of Bessel
processes (BESQ) to describe the evolution of masses. We are going
to carry it
 out in this section.

For $x\geq 0$ and $\delta\geq 0$ the {\it square of
$\delta$-dimensional Bessel process} starting at $x$, denoted by
$\text{BESQ}^\delta(x)$, is a non-negative valued process $\xi$
which solves the following stochastic differential equation
\[\xi_t=x+2\int_0^t\sqrt{\xi_s}dB_s+\delta t, \]
where $B$ is a one-dimensional Brownian motion. The Laplace
transform for $\xi$ is given by
\begin{equation}\label{lap-besq}
\bP\left[\exp\{-\lambda \xi_t\}\right]=\frac{1}{(1+2\lambda
t)^{\frac{\delta}{2}}} \exp\left\{-\frac{\lambda x}{1+2\lambda
t}\right\}.
\end{equation}
Notice that the Feller's branching process is just a
$\text{BESQ}^0$. We refer to Chapter XI in \cite{ReYo91} for a more
detailed introduction on the Bessel processes.

It is easy to see from (\ref{lap-besq}) that $\text{BESQ}^\delta(x)$
is additive in both $\delta$ and $x$; i.e. if $\{\xi_i,i=1,\ldots,m
\}$ is a sequence of independent processes such that each $\xi_i$ is
a $\text{BESQ}^{\delta_i}(x_i)$. Then $\sum_{i=1}^m\xi_i$ is a
$\text{BESQ}^{\sum_{i=1}^m\delta_i}(\sum_{i=1}^m x_i)$.

Now we are going to modify the process $Z$ defined in Section 3 by
letting the masses of the particles be governed by the BESQ
processes. Since the dimension is an additional parameter for BESQ,
we need to introduce another measure-valued process $\Delta$ to
describe the evolution of the dimension.

As in Section 3, we first consider two systems of interacting
particles. Given a finite measure $Z_0$ on $\bR$, for any $m$,
choose $x_1,\ldots,x_m$ to be i.i.d. random variables with a common
distribution $\bar Z_0:=Z_0/Z_0(\bR)$. Let $(X_1,\ldots,X_m)$ be an
$m$-dimensional coalescing Brownian motion starting at
$(x_1,\ldots,x_m)$.

Given another finite measure $\De_0$  on $\bR^+$ with a finite
``mean'' $\mu:=\int_0^\infty \de\De_0(d\de)$, let
$\de_1,\ldots,\de_m$ be i.i.d. random variables with a common
distribution $\bar\De:=\De_0/\De_0(\bR)$. We further suppose that
and $(X_i)$ and $(\de_i)$ are independent.

Put $\bar{\de}:=\De_0(\bR)$ and $\bar{z}:=Z_0(\bR)$. Let
$(\xi_1,\ldots,\xi_m)$ be a collection of $m$ independent
$\text{BESQ}^{\de_i\bar{\de}/m}(\bar{z}/m)$ processes. Then
\[Z_t^{(m)}:=\sum_{i=1}^m \xi_i(t)\delta_{X_i(t)}\]
and
\[\De_t^{(m)}:=\frac{1}{m}\sum_{i=1}^m\de_i\bar{\de}\delta_{X_i(t)}\]
define two  $M_F(\bR)$-valued processes.

Similar to Lemma \ref{tight} we can show that both $\{Z^{(m)}\}$ and
$\{\De^{(m)}\}$ are C-relatively compact in $D(M_F(\bR))$.  They
have unique weak limits by Theorem \ref{super-bess}, which we will
prove shortly.

Let $Z$ and $\De$ be the weak limits for $\{Z^{(m)}\}$ and
$\{\De^{(m)}\}$. Intuitively, $\{(Z_0(B),\Delta_0(B)):
B\in\cB(\bR)\}$ describes the initial mass-dimension distribution on
$\bR$, and $\{(Z_t(B),\Delta_t(B)): B\in\cB(\bR), 0\leq t<\infty\}$
describes the simultaneous mass-dimension evolution for such a
model, which we call a {\it super square of Bessel process with
spatial coalescing Brownian motion}.

For any nonnegative constants $\al_j, \be_j, j=1,\ldots, n$ and
$t>0$, put
\[I_t(x):=\sum_{j=1}^n\al_j1\{]Y_{2j-1}(t),Y_{2j}(t)]\}(x)\] and
\[J_t(x):=\sum_{j=1}^n\be_j1\{]Y_{2j-1}(t),Y_{2j}(t)]\}(x),\] where,
as usual, $(Y_1,\ldots,Y_{2n})$ is an $2n$-dimensional coalescing
Brownian motion starting at $(y_1,\ldots,y_{2n})$. The next result
determines the joint distribution for $(Z_t(B),\Delta_t(B)),
B\in\cB(\bR)$.

\begin{Theo}\label{super-bess}
For any $\al_j\geq 0, \be_j\geq 0, j=1,\ldots,n$ and any
$y_1\leq\ldots\leq y_{2n}$, we have
\begin{equation}\label{gen-lap}
\begin{split}
&\bP\left[\exp\left\{-\sum_{j=1}^n\al_j Z_t(]y_{2j-1},y_{2j}])-\sum_{j=1}^n
\be_j\De_t(]y_{2j-1},y_{2j}])\right\}\right]\\
&=\bP\left[\exp\left\{-\int_{-\infty}^\infty
Z_0(dx)\left(\frac{\mu}{2\bar{z}}\ln\left(1+2tI_t(x)\right)
+\frac{I_t(x)}{1+2 tI_t(x)}+\frac{\mu}{\bar{z}}
J_t(x)\right)\right\}\right].
\end{split}
\end{equation}
\end{Theo}

\begin{proof}

To prove (\ref{gen-lap}), again, we first fix $(\de_1,\ldots,\de_m)$
and $(\xi_1,\ldots,\xi_m)$. It follows that
\begin{equation}\label{lap2}
\begin{split}
&\bP[\exp\{-\sum_{j=1}^n \al_j Z^{(m)}_t(]y_{2j-1},y_{2j}])
-\sum_{j=1}^n \be_j\De^{(m)}_t(]y_{2j-1},y_{2j}])\}]\\
&=\bP\left[\exp\left\{-\sum_{j=1}^n \left(\sum_{i=1}^m
\al_j \xi_i(t)1\{]y_{2j-1},y_{2j}]\}(X_i(t))+
\sum_{i=1}^m \frac{\be_j\de_i\bar{\de}}{m}1\{]y_{2j-1},y_{2j}]\}(X_i(t))\right)\right\}\right]\\
&=\bP\left[\exp\left\{-\sum_{j=1}^n \left(\sum_{i=1}^m
\al_j\xi_i(t)1\{]Y_{2j-1}(t),Y_{2j}(t)]\}(x_i) \right.\right.\right.\\
&\qquad\qquad\qquad\qquad\qquad\left.\left.+\sum_{i=1}^m\frac{\be_j\de_i\bar{\de}}{m}1\{]Y_{2j-1}(t),Y_{2j}(t)]\}
(x_i)\left.\right)\right\}\right]\\
&=\bP\left[\exp\left\{-\sum_{i=1}^m\xi_i(t)I_t(x_i)- \sum_{i=1}^m
\frac{\de_i\bar{\de}}{m}J_t(x_i)\right\}\right].
\end{split}
\end{equation}
We then  fix $(x_i)$ and $(\de_i)$, and  take expectations with
respect to $(\xi_1,\ldots,\xi_m)$. By (\ref{lap-besq}) the right
hand side of (\ref{lap2}) is equal to
\begin{equation}\label{lap4}
\begin{split}
&\bP\left[\prod_{i=1}^m
\left(1+2tI_t(x_i)\right)^{-\frac{\delta_i\bar{\de}}{2m}}\exp\left\{-\frac{\bar{z}I_t(x_i)}{m(1+2
tI_t(x_i))}- \frac{\de_i\bar{\de} J_t(x_i)}{m}\right\} \right]\\
&\quad= \bP\left[\left(\int_0^\infty \bar
Z_0(dx)\int_0^\infty\bar\De_0(d\de)
\left(1+2tI_t(x)\right)^{-\frac{\delta\bar\de}{2m}}\exp\left\{-\frac{\bar{z}I_t(x)}{m(1+2
tI_t(x))}- \frac{\de \bar{\de} J_t(x)}{m}\right\}\right)^m \right].
\end{split}
\end{equation}

Let $m\goto\infty$ in (\ref{lap4}). We finally have
\begin{equation*}
\begin{split}
&\bP\left[\exp\left\{-\sum_{j=1}^n\al_j Z_t(]y_{2j-1},y_{2j}])-\sum_{j=1}^n
\be_j\De_t(]y_{2j-1},y_{2j}])\right\}\right]\\
&=\lim_{m\goto\infty}\bP\left[\exp\left\{-\sum_{j=1}^n\al_j Z_t^{(m)}(]y_{2j-1},y_{2j}])
-\sum_{j=1}^n\be_j\De_t^{(m)}(]y_{2j-1},y_{2j}]) \right\}\right]\\
&= \lim_{m\goto\infty}\bP\left[\left\{\int_{\bR\times\bR^+}\bar
Z_0(dx)\bar\De_0(d\de)\left(1-\frac{\delta\bar{\de}}{2m}\ln\left(1+2tI_t(x)\right)\right)
\left(1-\frac{\bar{z}I_t(x)}{m(1+2tI_t(x))}-\frac{\de\bar{\de} J_t(x)}{m}\right)\right\}^m\right]\\
&= \lim_{m\goto\infty}\bP\left[\left\{\int_{-\infty}^\infty \bar
Z_0(dx)\int_0^\infty\bar\De_0(d\de)\left(1-\frac{\delta\bar{\de}}{2m}\ln\left(1+2tI_t(x)\right)
-\frac{\bar{z}I_t(x)}{m(1+2tI_t(x))}-\frac{\de\bar{\de} J_t(x)}{m}\right)\right\}^m\right]\\
&= \bP\left[\exp\left\{-\int_{-\infty}^\infty \bar
Z_0(dx)\int_0^\infty\bar\De_0(d\de)\left(\frac{\delta\bar{\de}}{2}\ln\left(1+2tI_t(x)\right)
+\frac{\bar{z}I_t(x)}{1+2 tI_t(x)}+\de\bar{\de} J_t(x)\right)\right\}\right]\\
&=\bP\left[\exp\left\{-\int_{-\infty}^\infty
Z_0(dx)\left(\frac{\mu}{2\bar{z}}\ln\left(1+2tI_t(x)\right)
+\frac{I_t(x)}{1+2 tI_t(x)}+\frac{\mu}{\bar{z}}
J_t(x)\right)\right\}\right].
\end{split}
\end{equation*}

\end{proof}

\begin{Rem}
Notice that $\Delta$ is just the process $Z$ in Theorem
\ref{lap-fun} with $\gamma=0$.
\end{Rem}

The generalized model considered in this section will not die out
for $\mu>0$. Many of the properties in Section 3 and 4 can be
discussed in a similar fashion. But we leave the details to the
interested readers.

\noindent {\bf Acknowledgement}: The author is grateful to Steven
Evans for a suggestion that results in Section 5 of this paper.
The author also thanks Carl M${\Ddot{\text{u}}}$ller for a helpful
comment.

\bibliographystyle{alpha}



\end{document}